\documentclass[11pt,a4paper]{article}

\usepackage{mathtools}
\usepackage{amsthm}
\usepackage{graphicx}
\usepackage{subfloat}
\usepackage{filecontents}
\usepackage{lipsum}
\usepackage{amsmath,bm}
\usepackage{blindtext}
\usepackage{amsmath}
\usepackage{amssymb}
 \usepackage{paralist}
  \usepackage{graphics}
  \usepackage{filecontents}
   \usepackage{epsfig}
\usepackage{bbold}
    \usepackage[colorlinks=true]{hyperref}
 \hypersetup{urlcolor=blue, citecolor=red}

\usepackage{lineno,hyperref}

\usepackage{amsmath,amsthm,accents}
\newcounter{enumerConta}
\usepackage{hyperref}
\usepackage{graphicx}
\usepackage{subfigure}
\usepackage{enumitem}
\usepackage{kbordermatrix}

\newtheorem{theorem}{Theorem}[section]
\newtheorem{corollary}{Corollary}

\newtheorem{lemma}[theorem]{Lemma}

\theoremstyle{definition}
\newtheorem{definition}[theorem]{Definition}
\newtheorem{remark}{Remark}

\newcommand{\ones} {\mathbb{1}}

\title{Spectra of hyperstars on public transportation networks}

\author{Eleonora Andreotti \thanks{Division of Vehicle Safety, Department of Mechanics and Maritime Sciences, Chalmers University of Technology, SE-412 96 Göteborg, Sweden} }

\begin{document}

\maketitle


\begin{abstract}
The purpose of this paper is to introduce a model to study structures which are widely present in public transportation networks. We show that, through hypergraphs, one can describe these structures and investigate the relation between their spectra. To this aim, we extend the structure of $(m,k)-$stars on graphs to hypergraphs: the $(m,k)$-hyperstars on hypergraphs. Also, by giving suitable conditions on the hyperedge weights we prove the existence of matrix eigenvalues of computable values and multiplicities, where the matrices considered are Laplacian, adjacency and transition matrices. By considering separately the case of generic hypergraphs and uniform hypergraphs, we prove that two kinds of vertex set reductions on hypergraphs with $(m,k)$-hyperstar are feasible, keeping the same eigenvalues with reduced multiplicity. Finally, some useful eigenvectors properties are derived up to a product with a suitable matrix, and we relate these results to Fiedler spectral partitioning on the hypergraph.

\end{abstract}

\section{Introduction}

In the present work we focus on structures which are typically present in public transportation networks.\\
Many real social, chemical and biological relations can be represented as hypergraphs \cite{2009PLSCB...5E0385K,Estradahypergraph,1303205,Mulas,Mulas2,Mulas3,Mulas4, KONSTANTINOVA2001365}. In fact, hypergraphs are a fundamental tool for studying objects that cannot be characterized by simple binary relationships. 
In \cite{Mulas}, for example, the authors focus on chemical reactions which involve multiple atoms simultaneously.\\
In order to study relationships among multiple objects, a representation that best describes the properties of the structures is fundamental, without losing information on the $M$-relations (where $M=2$ in simple graphs ) and at the same time without necessarily assigning rigid roles to the entities. For this purpose, hypergraphs turn out to be a thorough tool. \\
Specifically, for transportation networks, defining each public transport line through its stations allows us to have a complete picture of the service provided by the city.\\
\emph{How can we modify a public transportation network while leaving the service unchanged?} -- one has to keep this question in mind when aiming to either add a station, remove a station, change stations of a line or eliminate the line itself.\\

In this context, the Laplacian formalism, as well as its spectrum, can be used to find many useful properties of the hypergraph. In particular,
studying isospectral hypergraphs means maintaining some properties of the structure, such as the number of connected components, the bipartiteness, the size of the graph, etc. For more details, we refer the reader to \cite{CHAN2020416, orientedhyp2013,orientedhyp2014,orientedhyp2018,orientedhyp2019-2,orientedhyp2019-3, Mulas} .

In this framework, the aim of this work is to study spectral properties of hypergraphs. A special focus is given to uniform hypergraphs, that have a large use in more applicative areas, such as biology and social sciences \cite{kong2019a, Zhang_2010, SHEPHERD1990395}; and to hyperstars, that represent structures which are widely present in transportation networks.\\

Together with the spectral properties of hypergraphs we shall also extend some results on Fiedler's spectral partitioning, in particular we shall extend results obtained in the previous work to the case of hypergraphs \cite{Andreotti18}.\\

The paper is organized as follows: we begin by stating the used terminology and by giving some preliminary remarks (Section \ref{sec:1}), in Section \ref{sec:2} we extend the results obtained in \cite{Andreotti18} by generalizing the class of $(m,k)$-star on graphs to the $(m,k)$-hyperstar on hypergraphs. 
In Section \ref{sec:3} we define two reduced $(m, k)$-hyperstars in hypergraphs classes: in the first case the reduction consists in removing some vertices, but keeping the hyperedges (which will simply be reduced by the number of vertices removed), in the second case we remove some vertices together with the hyperedges that contain them. In both cases we show that it is possible to keep the same spectrum of the initial hypergraph. 
Finally, in Section \ref{sec:4} we draw some conclusions.

\section{Notations and preliminary remarks}\label{sec:1}

We consider an undirected weighted connected hypergraph $\mathcal H:=(\mathcal V, \mathcal E, w)$, where the $N$ vertices in $\mathcal V$ are joined by the $M$ hyperedges in $\mathcal E$, with weight function: $w:\mathcal E\rightarrow \mathbb R^{+}.$ Let the \textit{rank} and the \textit{anti-rank} of $\mathcal H$ be the maximum  and the minimum cardinality of the edges in the hypergraph, respectively.
If all hyperedges have the same cardinality $p$ (i.e. if the rank and the anti-rank of $\mathcal H$ are equal to $p$, the hypergraph is said to be $p$-\textit{uniform}, \cite{berge1985graphs}.\\
Exactly as for simple graphs, an hypergraph with $N$ vertices and $M$ hyperedges may be defined by the incidence matrix $(H_{ve})_{v\in \mathcal{V}, e\in \mathcal{E}}$, i.e. by the matrix  of dimension $N\times M$  in which the columns correspond to the hyperedges while the rows correspond to the vertices of the hypergraph, and where
$$H_{ve} = \begin{cases} w(e), & \mbox{if vertex $v$ is contained in edge $e$ } (v\in e) \\
0 & \mbox{otherwise }  \end{cases}$$
where $v\in\mathcal V$ and $e\in \mathcal E$.\\

The degree of the vertex $v$ and the degree of the hyperedge $e$ are calculated
respectively as
$$
d(v) := \sum_{e \in \mathcal E} H_{ve}=\sum_{e \in \mathcal E, v\in e} w(e)$$

$$d(e) := \sum_{v \in \mathcal V}H_{ve}=|e| w(e).$$

We define $D_v\in Sym_N(\mathbb R^+_0)$ and $D_e\in Sym_M(\mathbb R^+_0)$ as the two diagonal matrices such that each diagonal entry corresponds to the vertex degree and to the hyperedge degree, respectively.\\
The adjacency matrix $A$ of hypergraph $\mathcal H$ is defined as

$$A := H^{1/2}(H^T
)^{1/2} - D_v.$$
Therefore, we can define the standard hypergraph Laplacian and the normalized hypergraph Laplacian matrices for hypergraph as follows

$$L := D_A - A$$

$$\mathcal L := I - D_A^{-1/2}AD_A^{-1/2}$$
where
$$D_A=diag(H^{1/2}(H^T
)^{1/2}\ones)-D_v$$
and $\ones$ is the ones-vector, \cite{Zhou06learningwith}.\\

Whenever we refer to the $k$-th eigenvalue of a Laplacian matrix (standard or normalized), we refer to the $k$-th eigenvalue according to an increasing order.

Furthermore, we observe that by defining the transition matrix $T$ as $T:=D_A^{-1}A$ we can link the spectrum of $T$ and the spectrum of $\mathcal L$. 

First of all, we observe that $T$ is similar to $\tilde A:=D_A^{-1/2}AD_A^{-1/2}$ via the invertible matrix $D_A^{1/2}:$
$$D_A^{-1/2}\tilde AD_A^{1/2}=D_A^{-1/2}D_A^{-1/2}AD_A^{-1/2}D_A^{1/2}=D_A^{-1}A=T.$$
Therefore $\sigma(T)=\sigma(\tilde A)$, where $\sigma(\cdot)$ is the spectrum of the considered matrix, and it is easy to prove that the following statements are equivalent
\begin{enumerate}[label=\textbf{S.\arabic*},ref=S.\arabic*]
\item $v$ is an eigenvector of $\tilde A$ with eigenvalue $\lambda$ 
\item $v^TD_A^{1/2}$ is a left eigenvector of $T$ with the eigenvalue $\lambda$
\item \label{tildeAL3}$D_A^{-1/2}v$ is a right eigenvector of $T$ with eigenvalue $\lambda$
\setcounter{enumerConta}{\value{enumi}}
\end{enumerate}
Thus, linking the spectrum of $T$ and the spectrum of $\mathcal L$ is equivalent to linking the spectrum of $\tilde A$ and the spectrum of $\mathcal L$, and we can easily prove that the following statments are equivalent
\begin{enumerate}[label=\textbf{S.\arabic*},ref=S.\arabic*]
\item \label{tildeAL1} $v$ is an eigenvector of $\tilde A$ with eigenvalue $\lambda$ 
\setcounter{enumi}{\value{enumerConta}}
\item \label{tildeAL4} $v$ is an eigenvector of $\mathcal L$ with the eigenvalue $1-\lambda.$
\end{enumerate}

For the classical results on Laplacian matrices, one may refer to \cite{Chung97, Cohen_Havlin:2010, Newman:2010:NI:1809753, Anderson85, MERRIS1994143}. For results on Laplacian matrices associated to hypergraphs, reference can be made to the book by A. Bretto \cite{Bretto:2013:HTI:2500991}.\\

Regarding the spectral partitioning of hypergraphs we refer to Zhou et al. \cite{Zhou:2006:LHC:2976456.2976657},
who generalized the methodology of spectral partitioning on undirected graphs to hypergraphs. In particular, we recall the Fiedler partitioning as given from the entries' sign of the second eigenvector of  its Laplacian matrix \cite{fiedler73,fiedler75}.

\section{Eigenvalues multiplicity in hypergraph matrices}\label{sec:2}
In the present section, 
we define the $(m,k)$-hyperstar, which generalizes the $(m,k)$-star \cite{Andreotti18}, and which, in its turn, generalizes the star \cite{DBLP:books/daglib/0070576}.
Together with the definitions, we also extend results obtained on the $(m, k)$-star, in particular
we extend Theorem (3.1) in \cite{Andreotti18} for hypergraphs. By defining weighted $(m,k)$-stars from hypergraphs, namely a weighted $(m,k)$-hyperstars, we are able to generalize the results obtained on multiple eigenvalues of Laplacian matrices, transition and adjacency matrices also for hypergraphs. 
\subsection{$(m,k)$-hyperstar: eigenvalues multiplicity}
We recall that a $(m,k)$-star is a graph $\mathcal G=(\mathcal V, \mathcal E,w)$ whose vertex set $\mathcal V$ can be written as the disjoined union of two subsets $\mathcal V_1$ and $\mathcal V_2$ of cardinalities $m$ and $k$ respectively, such that the vertices in $\mathcal V_1$ have no connections among them, and each of these vertices is connected
with all the vertices in $\mathcal V_2$: i.e
$$\forall i\in \mathcal V_1,\forall j\in \mathcal V_2,\quad (i,j)\in \mathcal E$$
$$\forall i,j\in \mathcal V_1, \quad (i,j)\notin \mathcal E.$$

We denote a $(m,k)$-star graph with partitions of cardinatilty $|\mathcal V_1|=m$ and $|\mathcal V_2|=k$ by $S_{m,k}.$\\
We define $(m,k)$-hyperstar and generalized $(m,k)$-hyperstar as follows:
\begin{definition}[$(m,k)$-hyperstar: $HS_{m,k}$ ]
A $(m,k)$-hyperstar is a hypergraph $\mathcal H=(\mathcal V, \mathcal E,w)$ whose vertex set $\mathcal V$ can be written as the disjoined union of two subsets $\mathcal V_1$ and $\mathcal V_2$, $\mathcal V=\mathcal V_1\dot{\cup}\mathcal V_2$, of cardinalities $m$ and $k$ respectively, such that 
$\exists P\in\mathcal P(\mathcal V_2) $
    with 
    \begin{itemize}
        \item $\bigcup_{\tilde e\in P}\tilde e=\mathcal V_2$,
        \item $\mathcal E=\{e \mid e=v_1\cup \tilde e , \tilde e \in P, v_1\in\mathcal V_1\}$,
        \item $\ w(\tilde e\cup v_i)=w(\tilde e\cup v_j) , \forall \tilde e\in P, \ v_i,v_j\in\mathcal V_1.$
    \end{itemize}

By $HS_{m,k}$ we denote a $(m,k)$-hyperstar of subsets $\mathcal V_1$ and $\mathcal V_2$ of cardinalities $|\mathcal V_1|=m$ and $|\mathcal V_2|=k$.

\end{definition}

\begin{definition}[Generalized $(m,k)$-hyperstar: $GHS_{m,k}$ ]
A generalized $(m,k)$-hyperstar is a hypergraph $\mathcal H=(\mathcal V, \mathcal E,w)$ whose vertex set $\mathcal V$ can be written as the disjoined union of two subsets $\mathcal V_1$ and $\mathcal V_2$, $\mathcal V=\mathcal V_1\dot{\cup}\mathcal V_2$, of cardinalities $m$ and $k$ respectively, and
$\forall v\in \mathcal V_1\exists P_{v}\in\mathcal P(\mathcal V_2) $
    such that 
    \begin{itemize}
        \item $\bigcup_{\tilde e\in P_{v}}\tilde e=\mathcal V_2$,
        \item $\mathcal E=\bigcup_{v\in\mathcal V_1}\{e \mid e=v\cup \tilde e , \tilde e \in P_v\}$,
        \item $\forall u\in\mathcal V_2$, $\ \sum_{\substack{u\in \tilde e,\\ \tilde e\in P_{v_i}}}w(\tilde e\cup v_i)=\sum_{\substack{u\in \tilde e,\\ \tilde e\in P_{v_j}}}w(\tilde e\cup v_j) , \forall v_i,v_j\in\mathcal V_1.$
    \end{itemize}

By $GHS_{m,k}$ we denote a generalized $(m,k)$-hyperstar of subsets $\mathcal V_1$ and $\mathcal V_2$ of cardinalities $|\mathcal V_1|=m$ and $|\mathcal V_2|=k$.
\end{definition}

\begin{remark}
A $(m, k)$- hyperstar is, trivially, a generalized $(m, k)$ -hyperstar such that $P_{v_i}=P_{v_j}, \forall v_i,v_j\in\mathcal V_1.$ Therefore, we shall consider generalized $(m, k)$-hyperstars to prove the results.
\end{remark}
Throughout this paper, we shall consider generalized $(m, k)$-hyperstars with $m,k\in \mathbb N$. When not else specified, we shall denote $GHS_{m,k}$ simply by $GHS$.\\

We define a \textit{generalized $(m,k)$-hyperstar on a hypergraph} $\mathcal H=(\mathcal V, \mathcal E,w)$ as the generalized $(m,k)$-hyperstar of partitions $\mathcal V_1$, $\mathcal V_2\subset \mathcal V$ such that only the vertices in $\mathcal V_2$ can be joined to the rest of the hypergraph $\mathcal V\setminus(\mathcal V_1\cup\mathcal V_2)$: i.e.\\

\begin{enumerate}[label=\textbf{(C.\arabic*)},ref=(C.\arabic*)]
\item \label{starcond2}
$\forall v\in \mathcal V_1\exists P_{v}\in\mathcal P(\mathcal V_2) $ and $\bar{\mathcal E}\subset\mathcal E$
    such that 
    \begin{itemize}
        \item $\bigcup_{\tilde e\in P_{v}}\tilde e=\mathcal V_2$,
        \item $\bar{\mathcal E}=\bigcup_{v\in\mathcal V_1}\{e \mid e=v\cup \tilde e , \tilde e \in P_v\}$,
        \item $\forall u\in\mathcal V_2$, $\ \sum_{\substack{u\in \tilde e,\\ \tilde e\in P_{v_i}}}w(\tilde e\cup v_i)=\sum_{\substack{u\in \tilde e,\\ \tilde e\in P_{v_j}}}w(\tilde e\cup v_j) , \forall v_i,v_j\in\mathcal V_1.$
    \end{itemize}
\item \label{starcond1} $\forall v_1\in \mathcal V_1, \forall v_2\in \mathcal V, v_1\neq v_2 , \quad \nexists e \in\mathcal E\setminus\bar{\mathcal E} \quad \mbox{ such that }\{v_1,v_2\}\subseteq e.$

\end{enumerate}

If there exists a generalized $(m,k)$-hyperstar on the hypergraph $\mathcal H$, then we say that the hypergraph $\mathcal  H$ has a generalized $(m,k)$-hyperstar.\\

\begin{figure}[!!h]
\begin{subfigure}{}
\includegraphics[width=10cm]{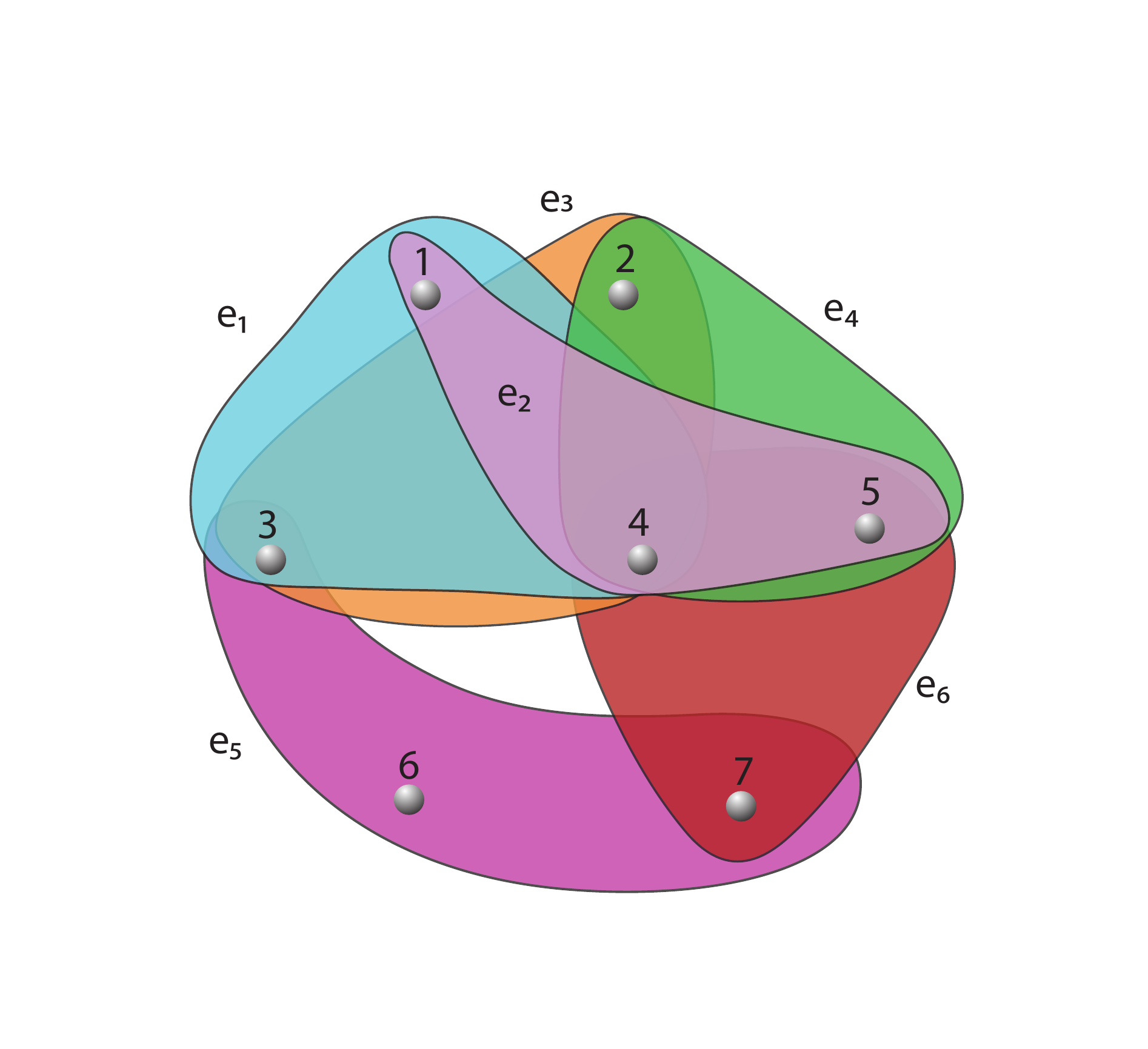}\vspace{-1cm}\end{subfigure}

\renewcommand{\kbldelim}{(}
\renewcommand{\kbrdelim}{)}
\[
  \text{H} = \kbordermatrix{
& e_1& e_2& e_3& e_4& e_5& e_6\\
1&1 & 2 & 0&0&0&0\\
2&0 & 0 & 1&2&0&0\\
3&1 & 0 & 1&0&1&0\\
4&1 & 2 & 1&2&0&3\\
5&0 & 2 & 0&2&0&3\\
6&0 & 0 & 0&0&1&0\\
7&0 & 0 & 0&0&1&3}
\]
\caption{A $HS_{2,3}$ on an hypergraph $\mathcal H=(\mathcal V, \mathcal E,w)$ and its incidence matrix. In this example $N=7$ and $M=6,\ \mathcal V_1=\{1,2\}$ and $\mathcal V_2=\{3,4,5\}$. The degree and weight of the $HS_{2,3}$ are $deg(HS_{2,3})=1$ and $w(HS_{3,4})=w_3+w_4+w_5=1+3+2=6$ respectively.}\label{fig:hyperstarunif}
\end{figure}
\begin{figure}[!!h]
\begin{subfigure}{}
\includegraphics[width=10cm]{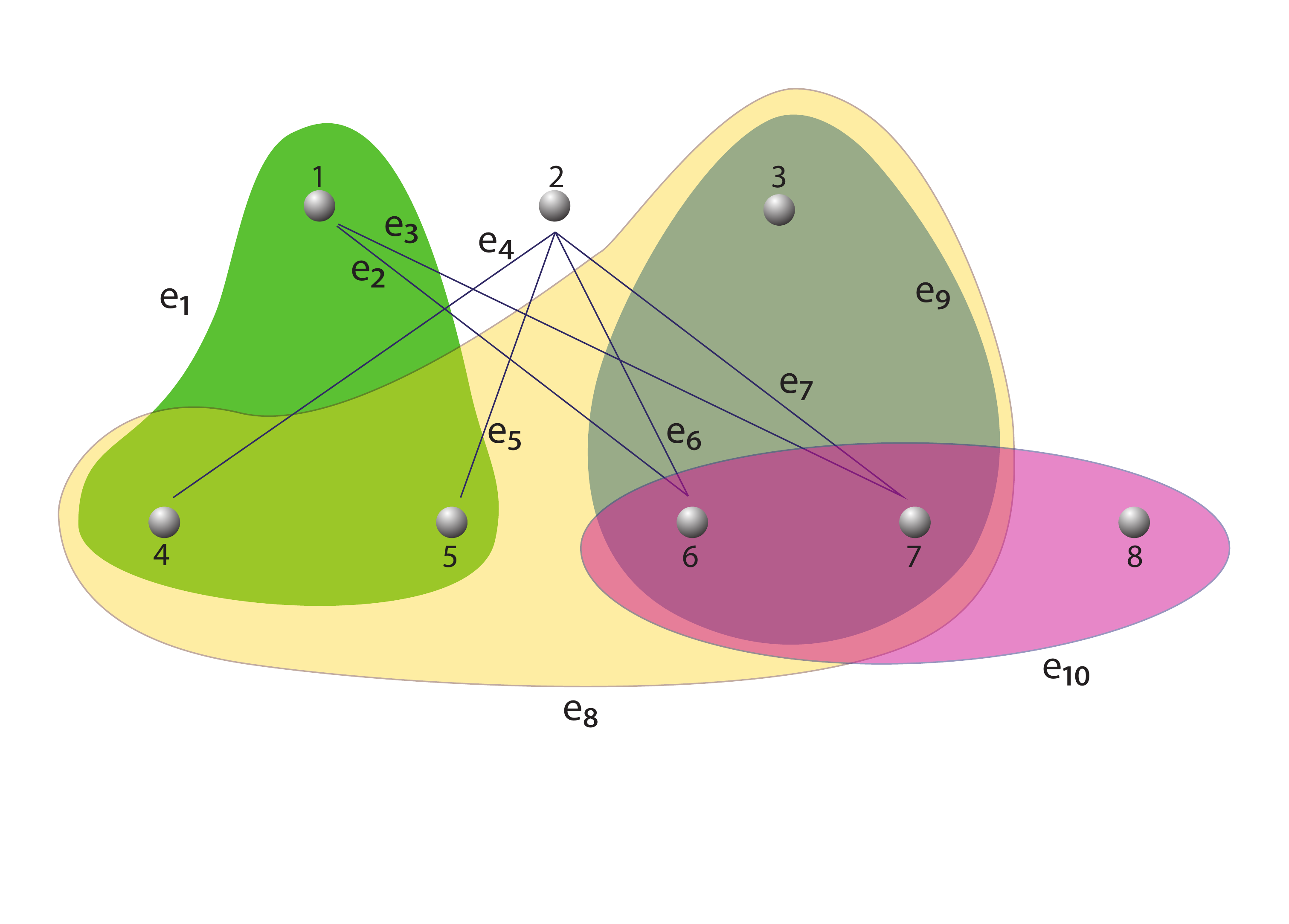}\vspace{-1cm}\end{subfigure}

\renewcommand{\kbldelim}{(}
\renewcommand{\kbrdelim}{)}
\[
  \text{H} = \kbordermatrix{
& e_1& e_2& e_3& e_4& e_5& e_6& e_7& e_8& e_9& e_{10}\\
1&1 & 2 & 2&0&0&0&0&0&0&0\\
2&0 & 0 & 0&1&1&2&2&0&0&0\\
3&0 & 0 & 0&0&0&0&0&1&1&0\\
4&1&0&0&1&0&0&0&1&0&0\\
5&1&0&0&0&1&0&0&1&0&0\\
6&0&2&0&0&0&2&0&1&1&3\\
7&0&0&2&0&0&0&2&1&1&3\\
8&0&0&0&0&0&0&0&0&0&3
}
\]
\caption{A $GHS_{3,4}$ on an hypergraph $\mathcal H=(\mathcal V, \mathcal E,w)$ and its incidence matrix. In this example $N=8$ and $M=10,\ \mathcal V_1=\{1,2,3\}$ and $\mathcal V_2=\{4,5,6,7\}$. The degree and weight of the $HS_{3,4}$ are $deg(GHS_{3,4})=2$ and $w(GHS_{3,4})=w_4+w_5+w_6+w_7=1+1+2+2=6$ respectively.}\label{fig:hyperstar2}
\end{figure}

In Figure \ref{fig:hyperstarunif} and Figure \ref{fig:hyperstar2} are represented an $(m, k)$-hyperstar and a generalized $(m, k)$-hyperstar on hypergraphs, respectively.\\

By defining the concepts of degree and weight of a generalized $(m,k)$-hyperstar we simplify the statement of the theorems on eigenvalues multiplicity.

\begin{definition}[Degree of a generalized $(m,k)$-hyperstar: $deg(GHS_{m,k})$]
The \textit{degree} of a generalized $(m,k)$- hyperstar is defined as follows $$deg(GHS_{m,k}):=m-1.$$
The degree of a set $\mathcal S$ of some $GHS$ such that $|\mathcal S|=l,$ is defined as the sum over each generalized $(m_i,k_i)$-hyperstar degree, $i\in\{1,...,l\}$, i.e.
$$deg(\mathcal S):=\sum_{i=1}^l deg(GHS_{m_i,k_i}).$$
\end{definition}

\begin{definition}[Weight of a generalized $(m,k)$-hyperstar: $w(GHS_{m,k})$]
The \textit{weight} of a generalized $(m,k)$-hyperstar {with} vertex set $\mathcal V_1\cup\mathcal V_2$, edge set $\mathcal E$ and weight function $w$, is defined as follows:
 $$w(GHS_{m,k}):=\sum_{v_2\in\mathcal V_2, \{v_1,v_2\}\subset e\in\mathcal E}w(e)\quad \quad \mbox{ for any } v_1\in\mathcal V_1.$$
\end{definition}
Before stating the extension to generalized $(m,k)-hyperstars$ of \cite[ Theorem 3.1]{Andreotti18}, we shall prove two useful Lemmas.
Given an hypergraph $\mathcal H=(\mathcal V,\mathcal E,w)$ associated with the adjacency matrix $A$, denoting $m_A(\lambda)$ the algebraic multiplicity of the eigenvalue $\lambda$ in $A$, the following Lemma holds.\\

\begin{lemma}\label{lemma:one}
Let $GHS$ be a generalized $(m,k)$-hyperstar of weight $w(GHS)$,
then
$$\exists \lambda \mbox{ such that } \lambda=0 \mbox{ and } m_{A}(\lambda)\geq deg(GHS).$$
\end{lemma}
\begin{proof}
Without loss of generality we consider only connected hypergraphs. 
{In fact,} if an hypergraph is not connected the same result holds, since the generalized $(m,k)$-hyperstar degree on the hypergraph is the sum of the generalized hyperstar degrees of the connected components and the characteristic polynomial of $A$ is the product of the characteristic polynomials of the connected components. \\
Under a suitable permutation of the rows and columns of {the} weighted incidence matrix $H$, we can label the vertices in $\mathcal V_1$ with the indices $1,...,m$, and the vertices in $\mathcal V_2$ with the indices $m+1,...,m+k$.\\
Let $v,u\in\mathcal V$, $v\neq u$, then the entry $v,u$ of the adjacency matrix is
$$A_{vu}=A_{uv}=\sum_{e\in\mathcal E}\sqrt {H_{ve} H_{ue}}=\sum_{e\in\mathcal E,\{v,u\}\subseteq e}w(e).$$

From condition \ref{starcond1}, if $v\in\mathcal V_1$ and $u\in\mathcal V_2$ (or $u\in\mathcal V_1$ and $v\in\mathcal V_2$), then $A_{vu}=A_{uv}=\sum_{\{u,v\}\subset e\in\mathcal E}=:w_u$ .\\
From condition \ref{starcond2}, if $v\in\mathcal V_1$ and $u\in\mathcal V\backslash\mathcal V_2$ (or $u\in\mathcal V_1$ and $v\in\mathcal V\backslash\mathcal V_2$), then $A_{vu}=A_{uv}=0$.\\

Let $v_1(A),...,v_m(A)$ be the rows corresponding to vertices in $\mathcal V_1=\{1,...,m\}$, then the adjacency matrix has the following form

\renewcommand{\kbldelim}{(}
\renewcommand{\kbrdelim}{)}
\[
  \text{A} = \kbordermatrix{
& v_1(A)&v_2(A)& & v_m(A)\\
v_1(A)& 0 & 0 &... & 0 &  w_{v_{m+1}}  &...& w_{v_{m+k}} & 0&...&0\\
v_2(A)&0 & 0 & \ddots & \vdots & \vdots  &...& \vdots & 0&...&0\\
& \vdots &\ddots & \ddots & 0 &  \vdots & ...& \vdots & 0&...&0\\
v_m(A)& 0  & ...& 0 & 0 &  w_{v_{m+1}}  &...& w_{v_{m+k}}& 0&...&0\\
&w_{v_{m+1}}& ...& ...  &  w_{v_{m+1}}  &  &  & & & &\\
 & \vdots & ... & ... & \vdots &  &  & & & &\\
& w_{v_{m+k}}& ... & ... &  w_{v_{m+k}} &  &  & & & & \\
& 0& ... & ... & 0 &  &  & & & & \\
& \vdots& ... & ... & \vdots &  &  A_{22} & & & \\
& 0& ... & ...&  0 &  &  & & & &  
}
\]
where the block $A_{22}$ is any $(N-m)\times(N-m)$ symmetric matrix with zero diagonal and nonnegative elements.\\ 
Because the matrix $ A$ has $m$ rows (and $m$ columns) $v_1( A),...,v_m( A)$ that are linearly dependent and such that $v_1( A)=...=v_m( A)$, then $m-1$ of these row vectors belong to the kernel of $A$.

Hence
$$\exists \mu_1,...,\mu_{m-1}\mbox{ eigenvalues of $A$ such that }\quad \mu_1=...=\mu_{m-1}=0.$$

\end{proof}

Similarly, let $L$ be the Laplacian matrix associated with the hypergraph $\mathcal H$. Denoting by $m_L(\lambda)$ the algebraic multiplicity of the eigenvalue $\lambda$ in $L$, the following Lemma holds.\\

\begin{lemma}\label{lemma:two}
Let $GHS$ be a generalized $(m,k)$-hyperstar of weight $w(GHS)$,
then
$$\exists \lambda \mbox{ such that } \lambda=w(GHS) \mbox{ and } m_{L}(\lambda)\geq deg(GHS).$$
\end{lemma}

\begin{proof}
Under a suitable permutation of the rows and columns of the weighted incidence matrix $H$, we can label the vertices in $\mathcal V_1$ with the indices $1,...,m$, and the vertices in $\mathcal V2$ with the indices $m+ 1,...,m+k$.
By Lemma \ref{lemma:one}, in the matrix $(-L+w(HS) I_N)$ there are the linearly dependent vectors $v_i, \ i\in\{1,...,m\}$, hence $m-1$ of these row vectors belong to $ ker(L-w(GHS) I_N)$ and
$$\exists \mu_1,...,\mu_{m-1}\mbox{ eigenvalues of $L-w(GHS) I_N$ such that }\quad \mu_1=...=\mu_{m-1}=0.$$
Let $\mu_i$ be one of these eigenvalues, then
$$0=det((L-w(GHS) I_N)-\mu_i I_N)=det(L-(w(GHS) +\mu_i )I_N)$$
so that $\lambda:=w(GHS)$ is an eigenvalue of $L$ with multiplicity greater or equal to $deg(GHS)$.\\
\end{proof}

We are now ready to enunciate the Theorem which extends \cite[ Theorem 3.1]{Andreotti18} to hypergarphs.
\begin{theorem}\label{th:one}
Let
\begin{itemize}
\item  $r$ be the number of $GHS$ in $\mathcal H$ with different weights, $w^1,...,w^r$, i.e. $w^i\neq w^j$ for each $i\neq j,$ where $ i,j\in\{1,...,r\};$
\item $\mathcal S_{w^i}$ be the set defined as follows
$$\mathcal S_{w^i}:=\{GHS\in \mathcal H | w(GHS)=w^i\}, \ i\in\{1,...,r\};$$
\end{itemize}
then for any $i\in\{1,...,r\},$
$$\exists \lambda \mbox{ such that } \lambda=w^i \mbox{ and } m_{ L}(\lambda)\geq \sum_{GHS\in \mathcal S_{w^i}}deg(GHS).$$
\end{theorem}
\begin{proof}
By using the same arguments as in Lemma \ref{lemma:two}, we can trivially prove the Theorem.
\end{proof}

Some corollaries on the normalized Laplacian matrix $\mathcal L$ and transition matrix $T$ can be obtained by similar proofs.

\begin{corollary}\label{cor:one}
Let
\begin{itemize}
\item  $r$ be the number of $GHS$ with different weights, $w^1,...,w^r$, i.e. $w^i\neq w^j$ for each $i\neq j,$ where $ i,j\in\{1,...,r\};$
\item $\mathcal S_{w^i}$ be the set defined as follows
$$\mathcal S_{w^i}:=\{GHS\in \mathcal H | w(GHS)=w^i\}, \ i\in\{1,...,r\};$$
\end{itemize}
then for any $i\in\{1,...,r\},$
$$\exists \lambda \mbox{ such that } \lambda=1 \mbox{ and } m_{\mathcal L}(\lambda)\geq \sum_{GHS\in \mathcal S_{w^i}}deg(GHS).$$
\end{corollary}

\begin{corollary}\label{cor:two}
Let
\begin{itemize}
\item  $r$ be the number of $GHS$ with different weight, $w^1,...,w^r$, i.e. $w^i\neq w^j$ for each $i\neq j,$ where $ i,j\in\{1,...,r\};$
\item $\mathcal S_{w^i}$ be the set defined as follows
$$\mathcal S_{w^i}:=\{GHS\in \mathcal H | w(GHS)=w^i\}, \ i\in\{1,...,r\};$$
\end{itemize}
then for any $i\in\{1,...,r\},$
$$\exists \lambda \mbox{ such that } \lambda=0 \mbox{ and } m_{T}(\lambda)\geq \sum_{GHS\in \mathcal S_{w^i}}deg(GHS).$$
\end{corollary}


\section{Generalized $(m,k)$-hyperstar dimensional reduction}\label{sec:3}

According to the previous results, we have defined a class of hypergraphs whose
Laplacian matrices have an eigenvalues spectrum with known multiplicities and
values. Now, our aim is to simplify the study of such hypergraphs by collapsing these
vertices into a single vertex replacing the original hypergraph with a reduced hypergraph. For this purpose we have defined two ways of collapsing 
the vertices. In the case of simple graphs these two modes are equivalent.\\
In subsection \ref{sub:3q} we define the generalized $(m,k)$-hyperstar $q$-reduction: this reduction consists in removing some vertices and reducing the cardinality of the hyperedges that contain them. In the case when $\mathcal H$ is a $p$-uniform hypergraph, then it is not guaranteed that the $q$-reduced hypergraph $\mathcal H^q$ is a $p$-uniform hypergraph too. 
In subsection \ref{sub:3q*} we define the generalized $(m,k)$-hyperstar $q_*$-reduction: this reduction consists in removing some vertices together with the hyperedges that contain them. In the case when the hypergraph $\mathcal H$ is a $p$-uniform hypergraph, then the $q_*$-reduced hypergraph $\mathcal H^{q_*}$ is a $p$-uniform hypergraph too. \\
After defining these two reduction classes of hypergraphs we will derive a spectrum correspondence between reduced and initial hypergraphs.\\
\subsection{Generalized $(m,k)$-hyperstar $q$-reduction}\label{sub:3q}
\begin{definition}[Generalized $(m,k)$-hyperstar $q$-reduced: $GHS^q_{m,k}$] A generalized $q$-reduced $(m,k)$-hyperstar is a generalized $(m,k)$-hyperstar with of vertex sets $\mathcal V_1,\mathcal V_2$, such that $q$ of its vertices in $\mathcal V_1$ are removed, $q<m$.\\
In other words: let $\mathcal H$ be a generalized $(m,k)$-hyperstar $(\mathcal V_1,\mathcal V_2, \mathcal E, w)$. A $GHS^q_{m,k}$ is defined for any choice $\{v_1,...,v_q\}\subset \mathcal V_1$ as the hypergraph $$(\{\mathcal V_1\backslash \{v_1,...,v_q\},\mathcal V_2\}, \mathcal E^q,w_{\big|\mathcal E^q}),$$ where $\mathcal E^q:=\{e| e:=\tilde e\backslash\{v_1,...,v_q\}, \tilde e\in\mathcal E\}.$ \\
Hence, the order (of the matrix) and the degree of the $GHS_{m,k}^q$ are $m+k-q$ and $m-q-1$, respectively. 
\end{definition}

\begin{figure}[!!h]\label{fig:hyperred}
\begin{subfigure}{}
\includegraphics[width=6.1cm]{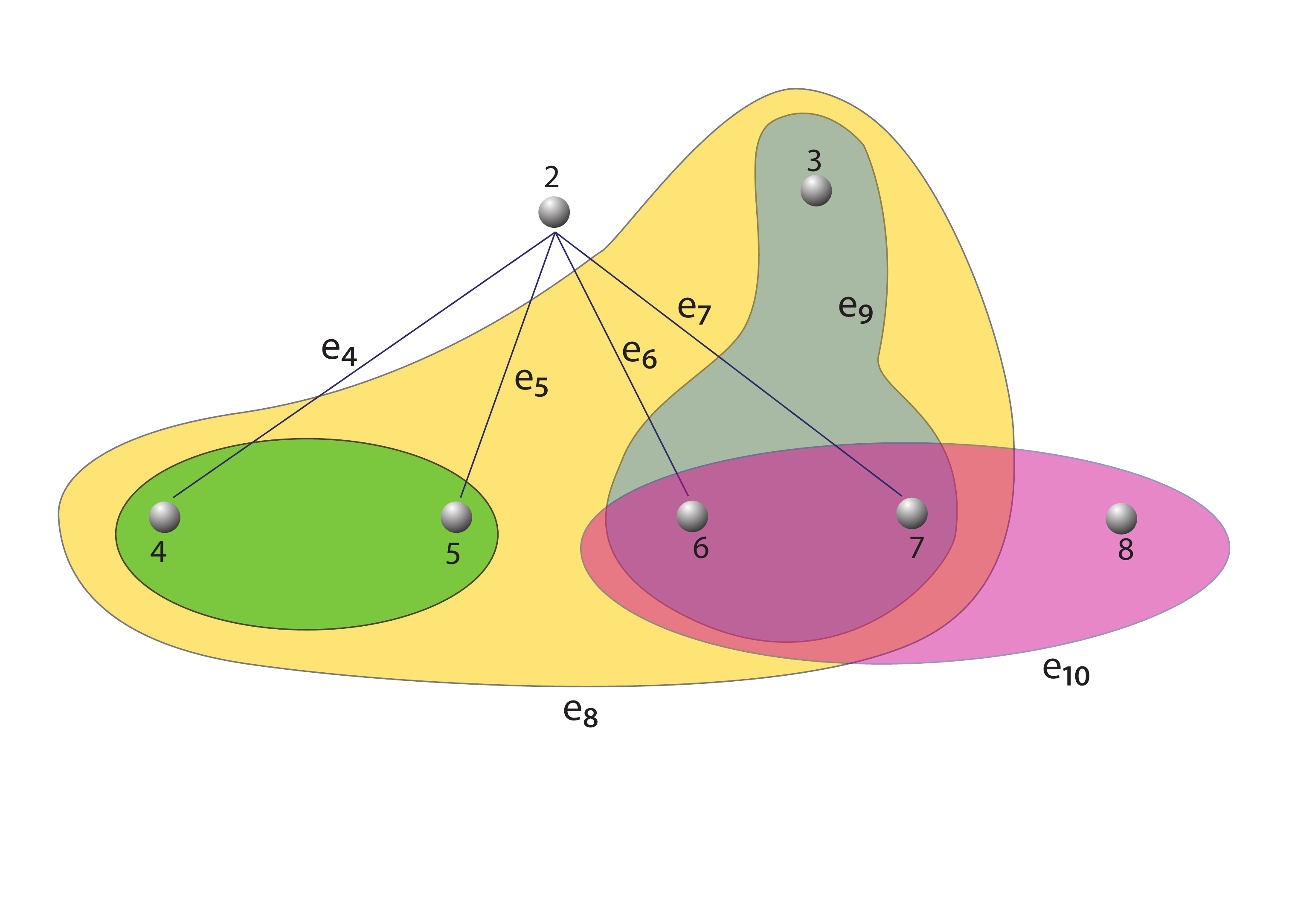}\vspace{-1cm}\end{subfigure}
\begin{subfigure}{}
\includegraphics[width=6.1cm]{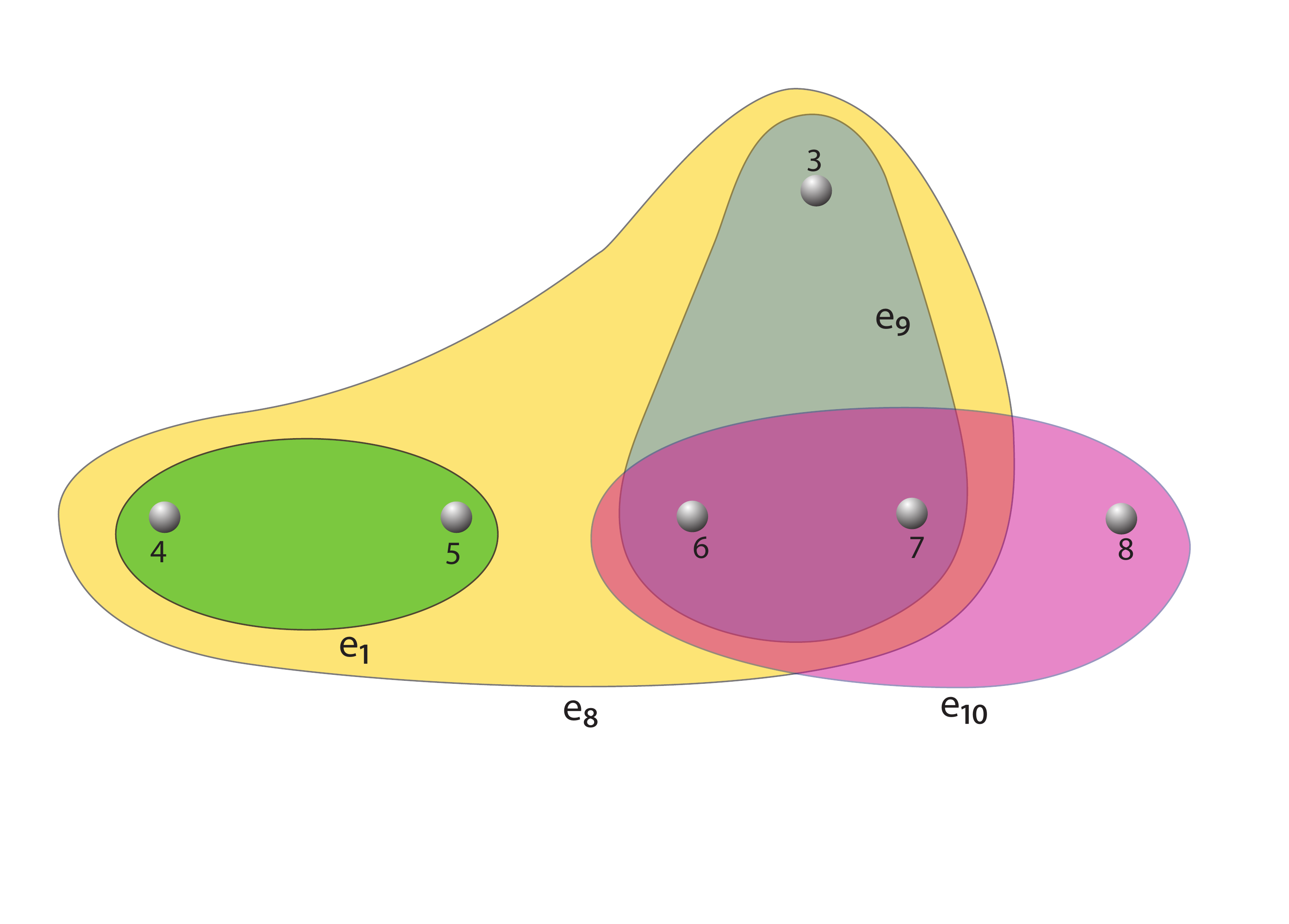}
\end{subfigure}
\caption{Examples of generalized $q$-reductions of the hypergraph $\mathcal H$ described in Figure \ref{fig:hyperstar2}: $GHS^1_{3,4}$ (left) and $GHS^2_{3,4}$ (right). }
\end{figure}

\begin{definition}[$q$-reduced hypergraph: $\mathcal H^q$] A $q$-reduced hypergraph $\mathcal H^q$ is obtained from an hypergraph $\mathcal H:=(\mathcal V,\mathcal E,w)$ with a generalized $(m,k)$-hyperstar (of vertex sets $\mathcal V_1,\mathcal V_2\subset\mathcal V$) by removing $q$ of the vertices in the set $\mathcal V_1$ of $\mathcal H$ and the set of hyperedges becomes $\mathcal E^q:=\{e| e:=\tilde e\backslash\{v_1,...,v_q\}, \tilde e\in\mathcal E\}$, where $\{v_1,...,v_q\}$ are the removed vertices.
Then 
$$\mathcal H^q:=(\mathcal V\backslash \{v_1,...,v_q\},\mathcal E^q,w_{\big|\mathcal E^q}).$$
\end{definition}

\begin{remark}
Whenever the hypergraph $\mathcal H$ is a $p$-uniform hypergraph, then it is not guaranteed that the $q$-reduced hypergraph $\mathcal H^q$ is a $p$-uniform hypergraph too. 
\end{remark}
 
Now we derive a spectrum correspondence between the hypergraphs $\mathcal H$ and $\mathcal H^q$.
\begin{definition}[Mass matrix of $GHS_{m,k}^q$] Let $\mathcal V_1,\mathcal V_2$ be the vertex sets of the hypergraph $GHS_{m,k}^q$, $q<m$. The \textit{mass matrix of} $GHS^q_{m,k},$ $\mathcal M$, is a diagonal matrix of order $m+k-q$ such that
$$\mathcal M_{vv} = \begin{cases} \frac{m}{m-q}, & \mbox{if } v\in\mathcal V_1\backslash\{v_1,...,v_q\}, \\
1 & \mbox{otherwise}
\end{cases}$$
\end{definition}

Similarly, we define the mass matrix $\mathcal M$ for an hypergraph $\mathcal H^q$, with a $GHS_{m,k}^q$, by means of a diagonal matrix of order $N-q$:

\begin{definition}[Mass matrix of $\mathcal H^q$] Let $\mathcal V$ be the vertex set of the hypergraph $\mathcal H$, $|\mathcal V|=N$, and $\mathcal V_1,\mathcal V_2$ be the vertex sets of the hypergraph $GHS_{m,k}^q$, $q<m$. The \textit{mass matrix of} $\mathcal H^q,$ $\mathcal M$, is a diagonal matrix of order $N-q$ such that
$$\mathcal M_{vv} = \begin{cases} \frac{m}{m-q}, & \mbox{if } v\in\mathcal V_1\backslash\{v_1,...,v_q\}, \\
1 & \mbox{otherwise}
\end{cases}$$
\end{definition}
For simplicity of notation we gave the definition of mass matrix of $\mathcal H^q$ with only one $GHS_{m,k}^q$,  but it can easily be extended to the case of multiple $GHS_{m,k}^q$.\\

\begin{theorem}[Generalized $(m,k)$-hyperstar adjacency matrix $q$-reduction theorem]\label{th:reduction1}
Let
\begin{itemize}
\item $\mathcal H$ be a hypergraph, on $N$ vertices, with a $GHS_{m,k}, \ m+q\leq N$,
\item $\mathcal H^q$ be the $q$-reduced hypergraph with a $GHS_{m,k}^q$ instead of $GHS_{m,k}$, on $N-q$ vertices,
\item $A$ be the adjacency matrix of $\mathcal H$,
\item $B$ be the adjacency matrix of $\mathcal H^q$,
\item $\mathcal M$ be the diagonal mass matrix of $\mathcal H^q$,
\end{itemize}
then
\begin{enumerate}
\item $\lambda$ is eigenvalue of $A \Leftrightarrow \ \lambda$ is eigenvalue of $\mathcal M B;$
\item There exists a matrix $K\in\mathbb R^{N\times (N-q)}$ such that $\mathcal M^{1/2}B\mathcal M^{1/2}=K^TAK$ and $K^TK=I$. Therefore, if $x$ is an eigenvector of $\mathcal M^{1/2}B\mathcal M^{1/2}$ for an
eigenvalue $\mu$, then Kx is an eigenvector of A for the same eigenvalue $\mu$.
\end{enumerate}
\end{theorem}

{ Before proving Theorem \ref{th:reduction1}, we recall a well known result for eigenvalues of
symmetric matrices, \cite{Hwang2004}.

\begin{lemma}[Interlacing theorem]
Let $A\in Sym_{N_A}(\mathbb R)$ with eigenvalues $\mu_1(A)\geq...\geq \mu_{N_A}(A).$ For $M<N$, let $K\in\mathbb R^{N_A,N_B}$ be a
matrix with orthonormal columns, $K^TK=I$, and consider the  $B=K^TAK$ matrix, with eigenvalues $\mu_1(B)\geq...\geq \mu_{N_B}(B).$
If
\begin{itemize}
\item the eigenvalues of $B$ interlace those of $A$, that is,
$$\mu_v(A)\geq\mu_v(B)\geq\mu_{N_A-N_B+v}(A), \quad v=1,...,N_B,$$
\item the interlacing is tight, that is, for some $0\leq u\leq N_B,$
$$\mu_v(A)=\mu_v(B), \ v=1,...,u\ \mbox{ and } \ \mu_v(B)=\mu_{N_A-N_B+v}(A), \ v=u+1,...,N_B$$
then $KB=AK.$
\end{itemize}
\end{lemma}

\begin{proof}
First we prove the existence of the $K$ matrix:\\
let $\mathcal P=\{P_1,...,P_{N-q}\}$ be a partition of the vertex set $\{1,...,N\}$.
 The \textit{characteristic matrix H} is defined as the matrix where the $u$-th column is the characteristic vector of $P_u$ ($u=1,...,N-q$).\\
Let A be partitioned according to $\mathcal P$
\[
A=\left(
\begin{array}{ccc}
A_{1,1} & \dots & A_{1,N-q} \\
\vdots &  & \vdots \\
A_{N-q,1} & \dots & A_{N-q,N-q}
\end{array}
\right),\]

where $A_{vu}$ denotes the block with rows in $P_v$ and columns in $P_u$.\\
The matrix $B=(b_{vu})$ whose entries $b_{vu}$ are the averages of the $A_{vu}$ rows, is called the \textit{quotient matrix} of $A$ with respect $\mathcal P$,
i.e. $b_{vu}$ denotes the average number of hyper-neighbours in $P_u$ of the vertices in $P_v$.\\
The partition is equitable if for each $v,u$, any vertex in $P_v$ has exactly $b_{vu}$ hyper-neighbours in $P_v$.
In such a case, the eigenvalues of the quotient matrix $B$ belong to the spectrum of $A$ ($\sigma(B)\subset\sigma(A)$) and the spectral radius
of $B$ equals the spectral radius of $A$: for more details cf. \cite{brouwer12}, chapter 2.\\
Also, we have the relations
$$\mathcal MB=H^TAH, \quad H^TH=\mathcal M.$$

Considering a $q$-reduced $(m,k)-$hyperstar with adjacency matrix $B$, we weight it by a diagonal mass matrix $\mathcal M$ whose diagonal entries
are all one except for the $m-q$ entries of the vertices in $\mathcal V_1$,

\begin{equation}\label{eq:diagM}
\mathcal M_{vv} = \begin{cases} \frac{m}{m-q}, & \mbox{if } v\in\mathcal V_1 \\ 1 & \mbox{otherwise } \end{cases},
\end{equation}

and we get
$$\mathcal MB\sim \mathcal M^{1/2}B\mathcal M^{1/2}=K^TAK, \quad K^TK=I,$$

where $K:=H\mathcal M^{1/2}.$
In addition to th.(\ref{th:one}), the eigenvalues of $\mathcal MB$ (with multiplicity) are also eigenvalues of $A$,
the adjacency matrix of the corresponding $HS_{m,k}$ hypergraph

$$
\sigma(\mathcal MB)\subset\sigma(A).
$$

Provided $q<m-1$, we get $\sigma(\mathcal MB)=\sigma(A)$, up to the multiplicity of the eigenvalue $\mu=0$.\\

Finally, if $x$ is an eigenvector of $\mathcal M^{1/2}B\mathcal M^{1/2}$ with eigenvalue $\mu$, then $Kx$ is an eigenvector of $A$ with the same eigenvalue $\mu$.\\

In fact from the equation
$$\mathcal M^{1/2}B\mathcal M^{1/2} x=\mu x,$$
taking into account that the partition is equitable, we have $K\mathcal M^{1/2}B\mathcal M^{1/2}=AK$ and
$$AKx=K\mathcal M^{1/2}B\mathcal M^{1/2}x=\mu Kx.$$

\end{proof}}

We obtain a similar result for the Laplacian matrix.

\begin{theorem}[Generalized $(m,k)$-hyperstar Laplacian matrix q-reduction theorem]\label{th:reduction2}
If
\begin{itemize}
\item $\mathcal H$ is an hypergraph, {on} N vertices, with a $GHS_{m,k}, \ m+q\leq N$,
\item $\mathcal H^q$ is the $q$-reduced hypergraph with a $GHS_{m,k}^q$ instead of $GHS_{m,k}$, of $N-q$ vertices,
\item $L(A)$ {is} the Laplacian matrix of $\mathcal H$,
\item $L(B)$ {is} the Laplacian matrix of $\mathcal H^q$,
\item $\mathcal M$ {is} the diagonal mass matrix of $\mathcal H^q$,
\end{itemize}
then
\begin{enumerate}
\item $\lambda$ is eigenvalue of $L(A) \Leftrightarrow \ \lambda$ is eigenvalue of $L(\mathcal M B);$
\item There exists a matrix $K\in\mathbb R^{N\times (N-q)}$ such that $\mathcal M^{1/2}B\mathcal M^{1/2}=K^TAK$ and $K^TK=I$.
Therefore, if $x$ is an eigenvector of $ \tilde L(\mathcal M B):=diag(\mathcal MB)-\mathcal M^{1/2}B\mathcal M^{1/2}$ for an eigenvalue $\lambda$, then $Kx$ is
an eigenvector of $L(A)$ for the same eigenvalue $\lambda$.
\end{enumerate}
\end{theorem}
The proof for the Laplacian version of the Reduction Theorem \ref{th:reduction1} is similar to that for the adjacency matrix,
in fact using the same arguments as in the proof of
\ref{th:reduction1}, we can say that 1. is true and that the $K$ matrix exists. {Hence} we prove directly only the second part of point 2. of the theorem.
\begin{proof}
Let $v$ be an eigenvector of $\tilde L(\mathcal M B):=diag(\mathcal MB)-\mathcal M^{1/2}B\mathcal M^{1/2}$ for an eigenvalue $\lambda$, then
$$ \tilde L(\mathcal M B)v=\lambda v.$$
Since $K\mathcal M^{1/2}B\mathcal M^{1/2}=AK$ and $diag(A)K=Kdiag(\mathcal MB)$, we obtain
$$L(A)Kx=diag(A)Kx-AKx=Kdiag(\mathcal MB)x-K\mathcal M^{1/2}B\mathcal M^{1/2} x=\lambda Kx.$$
\end{proof}

According to the previous results, an hypergraph with a generalized $(m,k)$-hyperstar and its $q$-reduced hypergraphs can be partitioned in the same way, up to removed vertices.\\

\begin{corollary}\label{cor:reduction}
Under the hypothesis of theorem \ref{th:reduction2}, if $x$ is a (left or right) eigenvector of $L(\mathcal MB)$ with eigenvalue $\lambda$,
then its entries have the same signs as the entries of the eigenvector $y$ of $L(A)$, with the same eigenvalue $\lambda$. \end{corollary}

\begin{proof} 
$ \tilde L(\mathcal MB)$ and $L(\mathcal MB)$ are similar by means of the matrix $\mathcal M^{1/2}$, in fact
\begin{eqnarray}
\mathcal M^{-1/2}L(\mathcal MB)\mathcal M^{1/2}&=&\mathcal M^{-1/2}diag(\mathcal MB)\mathcal M^{1/2}-\mathcal M^{-1/2}\mathcal MB\mathcal M^{1/2}\nonumber\\
&=&diag(\mathcal MB)-\mathcal M^{1/2}B\mathcal M^{1/2}\nonumber\\
&=&\tilde L(\mathcal MB). \nonumber
\end{eqnarray}

$L(\mathcal MB)$ preserves the sign of the eigenvectors of $\tilde L(\mathcal MB)$.\\
If $\tilde x$ an eigenvector of $\tilde L(\mathcal MB)$ of the eigenvalue $\lambda\in\sigma(\tilde L(\mathcal MB))$, then
\begin{eqnarray}
\tilde L(\mathcal MB) \tilde x=\lambda \tilde x & \Leftrightarrow & \mathcal M^{-1/2}L(\mathcal MB)\mathcal M^{1/2} \tilde x=\lambda \tilde x\nonumber\\
& \Leftrightarrow & L(\mathcal MB)\mathcal M^{1/2} \tilde x=\lambda \mathcal M^{1/2} \tilde x\nonumber
\end{eqnarray}
As a consequence, $x:=\mathcal M^{1/2} \tilde x$ is {an} eigenvector of $L(\mathcal MB)$ {for} the eigenvalue $\lambda,$ and $x_v=(\mathcal M\tilde x)_v$,
$$x_v=\sum_{r=1}^{N-q} \mathcal M_{vr}\tilde x_r=\mathcal M_{vv}\tilde x_v.$$

\end{proof}

From standard
results in linear algebra, we know that the
eigenvector $x$ of $L$ associated with its smallest nonzero eigenvalue is used to partition a graph as well as to partition a hypergraph \cite{Zhou:2006:LHC:2976456.2976657}. Hence, the vertex set
is partitioned into $V^+ = \{v \in \mathcal V |x(v) \geq 0\}$ and $V^- = \{v \in \mathcal V |x(v) < 0\}$. \\
In view of the previous result, we can partition the primary hypergraph $\mathcal H$ containing the $(m,k)$-hyperstar and the $q$-reduced hypergraph $\mathcal H^q$, weighted
by the matrix $\mathcal M$, in the same way except for the removed vertices.\\
As well as the Laplacian matrix (standard or normalized), other matrices are also used in spectral partitioning, for example in \cite{10.1007/978-3-662-45237-0_6} the authors show that the adjacency matrix can be more suitable for partitioning than other Laplacian matrices. Thanks to theorem \ref{th:reduction1}, one can apply the EVSA algorithm to the adjacency matrix of the primary hypergraph $\mathcal H$ containing the generalized $(m,k)$-hyperstar as well as using the adjacency matrix of the $q$-reduced hypergraph $\mathcal H^q$, linking the two results.\\

\subsection{$(m,k)$-hyperstar $q_*$-reduction}\label{sub:3q*}
In this section we focus on uniform hypergraphs and define a reduction that maintains the property of a uniform hypergraph.
In order to maintain the property of uniform hypergraph in the reduction, we give the following definitions

\begin{definition}[$p$-uniform $(m,k)$-hyperstar: $p$-$UHS_{m,k}$ ]
A $p$-uniform $(m,k)$-hyperstar is a hypergraph $\mathcal H=(\mathcal V, \mathcal E,w)$ whose vertex set $\mathcal V$ can be written as the disjoint union of two subsets $\mathcal V_1$ and $\mathcal V_2$, $\mathcal V=\mathcal V_1\dot{\cup}\mathcal V_2$, of cardinalities $m$ and $k$ respectively, and such that 
$\exists P\in\mathcal P(\mathcal V_2) $
    with 
    \begin{itemize}
        \item $\bigcup_{\tilde e\in P}\tilde e=\mathcal V_2$,
        \item $\forall \tilde e\in P, \mid \tilde e\mid=p-1$,
        \item $\mathcal E=\{e \mid e=v_1\cup \tilde e , \tilde e \in P, v_1\in\mathcal V_1\}$,
        \item $\ w(\tilde e\cup v_i)=w(\tilde e\cup v_j) , \forall \tilde e\in P, \ v_i,v_j\in\mathcal V_1.$
    \end{itemize}

By $p$-$UHS_{m,k}$ we denote a $p$-uniform $(m,k)$-hyperstar of subsets $\mathcal V_1$ and $\mathcal V_2$ of cardinalities $|\mathcal V_1|=m$ and $|\mathcal V_2|=k$.\\
When not else specified, we shall denote $p$-$UHS_{m,k}$ simply by $UHS_{m,k}$ or $UHS$.
\end{definition}

\begin{definition}[Uniform $(m,k)$-hyperstar $q_*$-reduced: $UHS^{q_*}_{m,k}$] A $q_*$-reduced uniform $(m,k)$-hyperstar is a uniform $(m,k)$-hyperstar of vertex sets $\mathcal V_1,\mathcal V_2$, such that $q$ of its vertices in $\mathcal V_1$ are removed together with all the hyperedges to which they belong.\\
In other words: let $\mathcal H$ be a $(m,k)$-hyperstar $(\mathcal V_1,\mathcal V_2, \mathcal E,w)$. A $UHS^{q_*}_{m,k}$ is defined for any choice $\{v_1,...,v_q\}\subset \mathcal V_1$ as the hypergraph $$(\{\mathcal V_1\backslash \{v_1,...,v_q\},\mathcal V_2\}, \mathcal E^{q_*},w_{\big|\mathcal E^{q_*}}),$$ where $\mathcal E^{q_*}:=\{e| v_i\notin e, i=\{1,...,q\}, e\in\mathcal E\}.$\\
The order and the degree of $UHS_{m,k}^{q_*}$ are $m+k-q$ and $m-q-1$, respectively.
\end{definition}
\begin{definition}[$q_*$-reduced hypergraph: $\mathcal H^{q_*}$] A $q_*$-reduced hypergraph $\mathcal H^{q_*}$ is obtained from a hypergraph $\mathcal H=(\mathcal V, \mathcal E, w)$ with a generalized $(m,k)$-hyperstar (of vertex sets $\mathcal V_1,\mathcal V_2\subset \mathcal V$) by removing $q$ of the vertices in the set $\mathcal V_1$ of $\mathcal H$ and the set of hyperedges becomes $\mathcal E^{q_*}:=\{e| v_i\notin e, i\in\{1,...,q\},  e\in\mathcal E\}$, where $\{v_1,...,v_q\}$ are the removed vertices. Then $\mathcal H^{q_*}:=(\mathcal V\backslash \{v_1,...,v_q\}, \mathcal E^{q_*},w_{\big|\mathcal E^{q_*}})$
\end{definition}

{We now derive a} spectrum correspondence between hypergraphs $\mathcal H$ and $\mathcal H^{q_*}$.
\begin{definition}[Vertices mass matrix of $UHS_{m,k}^{q_*}$] Let $\mathcal V_1, \mathcal V_2$ be the vertex sets of the hypergraph $UHS_{m,k}^q$, $q<m$. The \textit{vertices mass matrix of} $UHS_{m,k}^{q_*}$, $\mathcal M^*$, is a diagonal matrix of order $m+k-q$ such that
$$\mathcal M^*_{vv} = \begin{cases} \frac{m-q}{m}, & \mbox{if } v\in\mathcal V_1\backslash\{v_1,...,v_q\}, \\
1 & \mbox{otherwise}
\end{cases}$$
\end{definition}

\begin{definition}[Edges mass matrix of $UHS_{m,k}^{q_*}$] Let $\mathcal V_1,\mathcal V_2$ be the vertex sets of the hypergraph $UHS_{m,k}^{q_*}$, $q<m$.  The \textit{edges mass matrix of} $UHS_{m,k}^{q_*}$, $\mathcal N$, is a diagonal matrix of order $|\mathcal E^{q_*}|$ such that
$$\mathcal N_{ee} = \begin{cases} \frac{m}{m-q}, & \mbox{if }  e \cap \mathcal V_1\backslash\{v_1,...,v_q\}\neq \emptyset, \\
1 & \mbox{otherwise}
\end{cases}$$
\end{definition}

Similarly, we define the mass matrices $\mathcal M^*$ and $\mathcal N$ for an hypergraph $\mathcal H^{q_*}$, with one (or more) $UHS_{m,k}^{q_*}$. Even in this case, for simplicity of notation, we give the definition of mass matrices of $\mathcal H^{q_*}$ with only one $UHS_{m,k}^{q_*}$, but can easily be extended to the case of multiple $UHS_{m,k}^{q_*}$.\\

\begin{definition}[Vertices mass matrix of $\mathcal H^{q_*}$] Let $\mathcal V$ be the vertex set of the hypergraph $\mathcal H$, $|\mathcal V|=N$, and $\mathcal V_1,\mathcal V_2$ be the vertex sets of the hypergraph $UHS_{m,k}^{q_*}$, $q<m$. The \textit{vertices mass matrix of} $\mathcal H^{q_*},$ $\mathcal M^*$, is a diagonal matrix of order $N-q$ such that
$$\mathcal M^*_{vv} = \begin{cases} \frac{m}{m-q}, & \mbox{if } v\in\mathcal V_1\backslash\{v_1,...,v_q\}, \\
1 & \mbox{otherwise}
\end{cases}$$
\end{definition}

\begin{definition}[Edges mass matrix of $\mathcal H^{q_*}$] Let $\mathcal V$ be the vertex set of the hypergraph $\mathcal H$, $|\mathcal V|=N$, and $\mathcal V_1,\mathcal V_2$ be the vertex sets of the hypergraph $UHS_{m,k}^{q_*}$, $q<m$.  The \textit{edges mass matrix of} $\mathcal H^{q_*}$, $\mathcal N$, is a diagonal matrix of order $|\mathcal E^{q_*}|$ such that
$$\mathcal N_{ee} = \begin{cases} \frac{m}{m-q}, & \mbox{if }  e \cap \mathcal V_1\backslash\{v_1,...,v_q\}\neq \emptyset, \\
1 & \mbox{otherwise}
\end{cases}$$
\end{definition}








\begin{theorem}[Uniform $(m,k)$-hyperstar adjacency matrix $q_*$-reduction theorem]\label{th:reduction1Uni}
Let
\begin{itemize}
\item $\mathcal H$ be an hypergraph, on $N$ vertices, with a $UHS_{m,k}, \ m+q\leq N$, 
\item $\mathcal H^{q_*}$ be the $q_*$-reduced hypergraph with a $UHS_{m,k}^{q_*}$ instead of $UHS_{m,k}$, of $N-q$ vertices,
\item $A$ be the adjacency matrix of $\mathcal H$,
\item $I_{q_*}$ be the incidence matrix of $\mathcal H^{q_*}$,
\item $\mathcal M^*$ and $\mathcal N$ be the diagonal vertices and edges mass matrices of $\mathcal H^{q_*}$,
\end{itemize}
then
\begin{enumerate}
\item $\sigma(A)=\sigma(\mathcal M B), where B:=I_{q_*}^{1/2}\mathcal N(I_{q_*}^T)^{1/2}-diag(I_{q_*}^{1/2}\mathcal N(I_{q_*}^T)^{1/2})$
\item There exists a matrix $K\in\mathbb R^{N\times (N-q)}$ such that $\mathcal M^{1/2}B\mathcal M^{1/2}=K^TAK$ and $K^TK=I$. Therefore, if $x$ is an eigenvector of $\mathcal M^{1/2}B\mathcal M^{1/2}$ for an
eigenvalue $\mu$, then Kx is an eigenvector of A for the same eigenvalue $\mu$.
\end{enumerate}
\end{theorem}
\begin{proof}
By using the same arguments as in the proof of \ref{th:reduction1},  we can say that 1. and 2.
\end{proof}
We obtain a similar result for the Laplacian matrix.

\begin{theorem}[Uniform $(m,k)$-hyperstar Laplacian matrix $q_*$-reduction theorem]\label{th:reduction2Uni}
If 
\begin{itemize}
\item $\mathcal H$ is an hypergraph, of $N$ vertices, with a $UHS_{m,k}, \ m+q\leq N$,
\item $\mathcal H^{q_*}$ is the $q_*$-reduced hypergraph with a $UHS_{m,k}^{q_*}$ instead of $UHS_{m,k}$, of $N-q$ vertices,
\item $L(A)$ is the Laplacian matrix of $\mathcal H$,
\item $I_{q_*}$ is the incidence matrix of $\mathcal H^{q_*}$,
\item $\mathcal M^*$ and $\mathcal N$ are the diagonal vertices and edges mass matrices of $\mathcal H^{q_*}$,
\end{itemize}
then
\begin{enumerate}
\item $\sigma(L(A))=\sigma(L(\mathcal MB))$
\item There exists a matrix $K\in\mathbb R^{N\times (N-q)}$ such that $\mathcal M^{1/2}B\mathcal M^{1/2}=K^TAK$ and $K^TK=I$.
Therefore, if $x$ is an eigenvector of $\tilde L(\mathcal M B):=diag(\mathcal MB)-\mathcal M^{1/2}B\mathcal M^{1/2}$ for an eigenvalue $\lambda$, then $Kx$ is
an eigenvector of L(A) for the same eigenvalue $\lambda$.
\end{enumerate}
\end{theorem}
The proof for the Uniform versions of the Reduction Theorems are similar to the General one, in fact by using the same arguments as in the proofs of
\ref{th:reduction1} and \ref{th:reduction2}, we can prove the theorem.\\
According to the previous results, hypergraphs with $(m,k)$-hyperstars and {$q$-reduced hypergraphs} can be partitioned in the same way, up to the removed vertices.\\

\begin{corollary}\label{cor:reductionUni}
Under the hypothesis of theorem \ref{th:reduction2Uni}, if $x$ is a (left or right) eigenvector of $L(\mathcal MB)$ with eigenvalue $\lambda$,
then its entries have the same signs of the entries of the eigenvector $y$ of $L(A)$ with the same eigenvalue $\lambda$. \end{corollary}

\section{Conclusions}\label{sec:4}
{In this work, we have considered the problem of reducing the vertex set of a hypergraph while preserving spectral properies.}
 In presenting a vertex set reduction for hypergraphs, we defined the $(m,k)$-hyperstar, which generalizes the $(m,k)$-star \cite{Andreotti18}, and which, in its turn, generalizes the star \cite{DBLP:books/daglib/0070576}. We also generalized results concerning the value and the multiplicity of adjacency and Laplacian matrix eigenvalues, as it was done in \cite{Andreotti18} and \cite{Grone1994TheLS}. Unlike graphs with $(m,k)$-stars, for hypergraphs with $(m,k)$-hyperstars it is possible to define two different vertex set reductions, which lead to two different results on the reduction of the hypergraph: one can be performed on all types of hypergraphs, the other can be performed only on uniform hypergraphs.\\
The hyperstars introduced in this paper, together with {the generalization of structures} already defined for graphs, allow to describe structures {that are} present in transportation networks and to analyze when these structures have invariant characteristics, such as the spectrum or the sign of the eigenvectors. Thanks to these results we therefore know how to reduce the number of peripheral stations with an appropriate increase in the service provided, represented by the new hyperedge weights in the reduced graph. Future developments of the model concern the study of oriented and bipartite hypergraphs, in order to involve different means of transport.

\section*{Acknowledgments}The author would like to thank Raffaella Mulas (Max Planck Institute of Leipzig, Germany) for the helpful comments and discussions.

\bibliographystyle{unsrt}
\bibliography{biblio_coal}

\end{document}